
\magnification=1200
\input amssym.def
\input amssym.tex
 \font\newrm =cmr10 at 24pt
\def\bul{\raise .9pt\hbox{\newrm .\kern-.105em } }

 \def\fr{\frak}
 \font\sevenrm=cmr7 at 7pt

 \baselineskip 13pt
 
 \def\h{\hbox{ }}

 \def\r{{\fr r}}

 \def\n{{\fr n}}

 \def\ss{{\fr s}}
 
 \def\b{{\fr b}}
 
 \def\hh{{\fr h}}

 \def\g{{\fr g}}

 \def\<{\le}
 \def\>{\ge}

 \def\s{{\h\subset\h}}
 
 \def\vs{\vskip }

 \def\mapright#1
  {\smash{\mathop
  {\longrightarrow}
  \limits^{#1}}}

 \def\kk#1{{\kern .4 em} #1}
 \def\vs{\vskip 1pc}
 
\input macro.sty

\hsize = 27pc
\vsize = 45pc
\overfullrule = 0pt

\font\sevenrm=cmr10 at 7pt

\hsize=31pc
\vsize=45pc

\rm

\centerline {\bf Center of $U(\n)$, Cascade of Orthogonal Roots}
\vskip 3pt
\centerline{\bf and a Construction of Lipsman--Wolf}

\vskip 2pc

\centerline{\bf Bertram Kostant}

\hfill {\it Dedicated to Joe, a special friend and valued colleague}

\vskip 1.5pc 
\noindent {\bf Abstract}. Let
$G$ be a complex simply-connected semisimple Lie group and let $\g=
 \hbox{\rm Lie}\,G$. Let $\g = \n_- +\hh + \n$ be a triangular decomposition of
$\g$. One readily has that  $\hbox{\rm Cent}\,U(\n)$ is isomorphic to the ring 
 $S(\n)^{\n}$ of symmetric invariants. Using the cascade ${\cal B}$ of strongly orthogonal
roots, some time ago we  proved (see [K]) that $S(\n)^{\n}$ is a polynomial
ring $\Bbb C[\xi_1,\ldots,\xi_m]$ where $m$ is the cardinality of ${\cal B}$. 
The authors in [LW] introduce a very nice representation-theoretic 
method for the construction of certain elements in
$S(\n)^{\n}$. A key lemma in [LW] is incorrect but the idea is in
fact valid. In our paper here we modify the construction so as to yield 
these elements in $S(\n)^{\n}$ and use the [LW] result to prove a theorem 
 of Tony Joseph.


\vskip 2pt

\noindent{\bf MSC (2010) codes}: representation theory, invariant theory.
\vskip 3pt

\noindent {\bf 1. Introduction}
\vskip 1pc {\bf 1.1.} Let $\g$ be a complex semisimple Lie algebra and let $$\g = \n_- + \hh + \n$$
be a fixed triangular decomposition of $\g$. Let $\Delta\s \hh^*$ be
the set of $\hh$ roots in $\g$. The Killing form
$(x,y)$ on $\g$, denoted by ${\cal K}$, induces a nonsingular
bilinear form $(\mu,\nu)$ on $\hh^*$. For each $\varphi\in \Delta$
let $e_{\varphi}\in \g$ be a corresponding root vector. The
root vectors can and will be chosen so that
$(e_{\varphi},e_{-\varphi}) = 1$ for all roots $\varphi$. 

If $\ss\s \g$ is any subspace stable under $ \hbox{\rm ad}\,\hh$ let
$$\Delta(\ss)=\{\varphi\in \Delta\mid e_{\varphi}\in \ss\}.$$ The set
$\Delta_+$ of positive roots is then chosen so that $\Delta_+ =
\Delta(\n)$, and one puts $\Delta_- = -\Delta_+$. If $\ss$ is a
Lie subalgebra, then
$S(\ss)$ and
$U(\ss)$ are respectively the symmetric and enveloping algebras of
$\ss$. Our concern here is with the case where
$\ss =\n$.

Let $\b = \hh + \n$ so that $\b$ is a Borel subalgebra of $\g$.
Let $G$ be a Lie group such that $ \hbox{\rm Lie}\,G = \g$ and let $H,N,B$
be Lie subgroups corresponding, respectively, to $\hh,\n,\b$.
Then $S(\n)$ is a $B$-module since $B=HN$ normalizes $N$. Let $m$
be the maximal number of strongly orthogonal roots. Then we proved the following
some time ago, generalizing a result of Dixmier (case where $\g$ is of
type $A_{\ell}$),

 \vs {\bf Theorem A.} {\it There exists $\xi_i\in
S(\n)^N, i=1,\ldots,m$, so that $$S(\n)^N = \Bbb
C[\xi_1,\ldots,\xi_m] $$ is a polynomial ring in $m$-generators. 
 Furthermore, $$S(\n)^N \cong  \hbox{\rm Cent}\,U(\n)$$ so
that one has a similar statement for $\hbox{\rm Cent}\,U(\n)$.}

 \vs  We will present an algebraic-geometric proof of a much
stronger statement than Theorem A and relate it to a 
representation-theoretic construction, due to Lipsman--Wolf, of certain
elements in
$S(\n)^N$. See [K], [LW]. A key tool is the cascade ${\cal B}=
\{\beta_1,\ldots,\beta_m\}$ of orthogonal roots which will now
be defined.

 \vs {\bf 1.2.} Let $\Pi\s \Delta_+$ be the set of
simple positive roots. For any $\varphi\in \Delta_+$ and $\alpha\in
\Pi$ there exists a nonnegative integer $n_{\alpha}(\varphi)$  such
that $$\varphi = \sum_{\alpha\in \Pi}\,n_{\alpha}(\varphi)\,\alpha.$$
Let $$\Pi(\varphi) = \{\alpha\in \Pi\mid n_{\alpha}(\varphi) >0\}.$$
Then $\Pi(\varphi)$ is a connected subset of $\Pi$ and hence defines
a simple Lie subalgebra $\g(\varphi)$ of $\g$. We will say that
$\varphi$ is locally high if $\varphi$ is the highest root of
$\g(\varphi)$. Obviously the highest roots of all the simple
components of $\g$ are locally high. 

\vs {\bf Remark 1.} If $\g$ is of type
$A_{\ell}$, but only in this case,  are all $\varphi\in \Delta_+$ 
locally high. 

\vs Let $\varphi\in\Delta_+$ be locally high and let
$$\Pi(\varphi)^o = \{\alpha\in \Pi(\varphi)\mid (\alpha,\varphi)=
0\}; $$  let $\g(\varphi)^o$ be the semisimple Lie algebra having
$\Pi(\varphi)^o$ as its set of simple roots. We will then say that
a root $\varphi'\in \Delta_+$ is an offspring of $\varphi$ if
$\varphi'$ is the highest root of a simple component of
$\g(\varphi)^o$. 

\vs {\bf Remark 2.} One notes that an offspring of
a locally high root $\varphi$ is again locally high and that it is
strongly orthogonal to $\varphi$.

\vs A sequence of positive roots $$C =
\{\beta_1',\ldots,
\beta_k'\}$$ will be called a cascade chain if $\beta_1'$ is a
highest root of a simple component of $\g$, and if $1<j\leq k$, then
$\beta_j'$ is an offspring of $\beta_{j-1}'$. Now let ${\cal B}$ be
the set of all positive roots $\beta$ which are members of some
cascade chain. Let $W$ be the Weyl of $(\hh,\g)$.

 \vs  {\bf Theorem 1.} {\it The cardinality of
${\cal B}$ is $m$ and $${\cal B} = \{\beta_1,\ldots,\beta_{m}\}$$ is
a maximal set of strongly orthogonal roots. Furthermore, if
$s_{\beta_i}$ is the $W$-reflection of $\hh$
corresponding to $\beta_i$, then the long element $w_o$ of $W$ may be
given by $$w_o = s_{\beta_1}\cdots s_{\beta_m}.\eqno (1.1)$$}

${\cal B}$ is the cascade of orthogonal roots.

 \vs {\bf 1.3.} One has the
vector space direct sum $$\g = \n_- \oplus \b.\eqno (1.2.)$$ Let $P:
\g \to \n$ be the projection defined by (1.2). Since $\b$ is the
${\cal K}$-orthogonal subspace to $\n$ in $\g$ we may identify
$\n_-$ with the dual space $\n^*$ to $\n$, so that for $v\in \n_-$
and $x\in \n$, one has $\langle v,x\rangle = (v,x)$. The  coadjoint
action of $N$ on $\n_-$ may then be given so that if $u\in N$,
then on $\n_-$ 
$$\hbox{\rm Coad}\,u= P\,\hbox{\rm Ad}\,u.\eqno (1.3.)$$
In fact, using (1.2) the coadjoint action of $N$ on $\n_-$ extends to an
action of $B$ on $\n_-$, so that if $b\in B$ and $v\in \n_-$, one has
$b\cdot v = P\,\hbox{\rm Ad}\,b(v)$. In addition we can regard
$S(\n)$ as the ring of polynomial functions on $\n_-$. Since $B$
normalizes $N$ the natural action of $N$ on $S(\n)$ extends to
an action of $B$ on $S(\n)$ where if $f\in S(\n),\, b\in B$, and
$v\in\n_-$, one has $$(b\cdot f)(v) = f(b^{-1}\cdot v).\eqno
(1.4)$$

\vskip .5pc
Recalling $m = \hbox{card}\,{\cal B}$, let $\r$ be the
commutative
$m$-dimensional subalgebra of $\n$ spanned by $e_{\beta}$ for
$\beta\in {\cal B}$ and let
$R \s N$ be the commutative unipotent subgroup corresponding to
$\r$. In the dual space let $\r_-\s \n_-$ be the span of
$e_{-\beta}$ for
$\beta\in {\cal B}$. For any $z\in \r_-,\,\beta\in {\cal B}$, let
$a_{\beta}(z)\in 
\Bbb C$
be defined so that $$z = \sum_{\beta\in
{\cal B}}\,a_{\beta}(z)\,e_{-\beta},\eqno (1.5)$$ and let
$$\r_-^{\times} = 
\{\tau\in
\r_-\mid a_{\beta}(\tau)\neq 0,\,\,\forall
\beta\in {\cal B}\}.$$ As an algebraic subvariety of $\n_-$
clearly 
$$\r_-^{\times} \cong (\Bbb C^{\times})^m.\eqno (1.6)$$

Also for any $z\in \n_-$ let $O_z$ be the $N$-coadjoint orbit
containing $z$. Let $N_z\s N$ be the coadjoint isotropy subgroup at
$z$ and let $\n_z = \hbox{Lie}\,N_z$. Since the action is algebraic, 
$N_z$ is connected and hence as $N$-spaces $$O_z\cong N/N_z.\eqno
(1.7)$$

\vskip .5pc
{\bf Theorem 2.} {\it Let $\tau\in
\r_-^{\times} $. Then (independent of $\tau$) $N_{\tau} = R$ so that
(1.7) becomes $$O_{\tau}\cong N/R.\eqno (1.8)$$ 

In particular
$$\hbox{\rm dim}\,O_{\tau} = \hbox{\rm dim}\,\n - m\eqno (1.9)$$ and
$O_{\tau}$ is a maximal dimensional coadjoint orbit of $N$.} 

\vs Now consider
the action of $B$ on $\n_-$. In particular  consider the action of
$H$ on $\n_-$. Obviously $$\r_-^{\times} \cong (\Bbb
C^{\times})^m,\eqno (1.10)$$ and furthermore $\r_-^{\times}$ is an orbit
of
$H$.  In addition $H$ permutes the maximal $N$-coadjoint orbits
$O_{\tau},\,\tau\in \r_-^{\times}$. More precisely,

 \vs {\bf Theorem 3.} {\it For any $a\in H$ and $\tau\in \r_-^{\times}$, one has
$$a\cdot O_{\tau} = O_{a\cdot \tau}.\eqno (1.11)$$}

 {\bf 1.4.} If
$V$ is an affine variety, $A(V)$ will denote its corresponding affine
ring of functions. Note that $S(\n) = A(\n_-)$. Let $Q(\n_-)$ be the
quotient field of $S(\n)$. 

\vs {\bf Theorem 4.} {\it There exists a
unique Zariski open nonemtpy orbit $X$ of $B$ on $\n_-$. In particular
$${\overline X} = \n_-.\eqno (1.12)$$ Furthermore $X$ is an affine
variety so that $$S(\n_-)\s A(X)\s Q(\n_-).\eqno (1.13)$$ Moreover
$\n_-^{\times}\s X$, and in fact one has a disjoint union $$X =
\sqcup_{\tau\in \r_-^{\times}}O_{\tau}\eqno (1.14)$$ so that all
$N$-coadjoint orbits in $X$  are maximal and isomorphic to
$N/R$.}

\vs Let $\Lambda\s \hh^*$ be the $H$-weight
lattice and let $\Lambda_{\hbox{\sevenrm ad}}\s \Lambda$ be the root lattice.
Let
$\Lambda_{\cal B}\s \Lambda_{\hbox{\sevenrm ad}}$ be the sublattice 
 generated by
the cascade ${\cal B}$. Since the elements of ${\cal B}$ are mutually
orthogonal note that $$\Lambda_{\cal B} = \oplus_{\beta\in
{\cal B}}\,\Bbb Z\ \beta\eqno (1.15)$$ is a free $\Bbb Z$-module of
rank $m$. 

If $M$ is an $H$-module, let $\Lambda(M)\s \Lambda$ be the set of
$H$-weights occurring in $M$. Note that if $M$ is a
$B$-module, then $M^N$ is still an $H$-module. Recalling the definition
of $\r_-^{\times}$ and (1.6), note that
$$\eqalign{\Lambda(A(\r_-^{\times})) &=
\Lambda_{\cal B}\cr \hbox{and each weight}\,\, & \hbox{occurs with
multiplicity 1}.\cr}\eqno (1.16)$$ 

We can now give more information
about $X$ and its affine ring $A(X)$. Define a
$B$ action on $\r_-^{\times}$ by extending the $H$-action so
that $N$ operates trivially. Next define a $B$-action on $N/R$,
extending the $N$-action by letting $H$ operate by conjugation,
noting that $H$ normalizes both $N$ and $R$. With these structures
and the original action on $X$, we have the following. 

\vs {\bf Theorem 5.} {\it One has a $B$-isomorphism
$$X\to N/R\times \r_-^{\times}$$ of affine varieties so
that as $B$-modules $$A(X) \cong A(N/R)\otimes A(\r_-^{\times}).\eqno
(1.17)$$ Furthermore, taking $N$-invariants, one has an $H$-module
isomorphism $$A(X)^N\cong A(\r_-^{\times})\eqno (1.18)$$ so that, by
(1.16), $$\Lambda(A(X)^N) = \Lambda_{\cal B}\eqno (1.19)$$ and each $H$-weight
occurs with multiplicity $1$.}

\vs Recalling (1.13) one has the
$N$-invariant inclusions $$S(\n)^N \s A(X)^N
\s Q(n_-)^N\eqno(1.20)$$ of $H$-modules so that $$\Lambda(S(\n)^N) \s
\Lambda(A(X)^N)\s \Lambda(Q(n_-)^N).\eqno (1.21)$$ But since $S(\n)$ is
a unique factorization domain, any $u\in Q(\n_-)$ may be uniquely
written, up to scalar multiplication as $$u = f/g\eqno (1.22)$$ where
$f$ and
$g$ are prime to one another. Furthermore, it is then immediate (since
$N$ is unipotent) that if $u$ is $N$-invariant, one has $f,g\in
S(\n)^N$. If, in addition, $u$ is an $H$-weight vector, the same is
true of $f$ and $g$ so that, using Theorem 5, one readily
concludes the following.

 \vs {\bf Theorem 6.} {\it Every $H$-weight in
$\Lambda(S(\n)^N)$ occurs with multiplicity $1$ in $S(n)^N$. In fact
$\Lambda(Q(\n_-) = \Lambda_{\cal B}$ and every weight $\gamma$ in
$\Lambda(Q(\n_-)$ occurs with multiplicity $1$ in $Q(\n_-)^N$ and is
of the form $$\gamma= \nu-\mu\eqno (1.23)$$ where $\mu,\nu \in
\Lambda(S(\n)^N)$. }

\vs For any $\gamma\in \Lambda_{\cal B}$ let
$\xi_{\gamma}\in Q_{\n_-}^N$ be the unique (up to scalar
multiplication) $H$-weight vector with weight $\gamma$. Thus if
$\gamma\in \Lambda_{\cal B}$, we may uniquely  write (up
to scalar multiplication $$\xi_{\gamma} =\xi_{\nu}/\xi_{\mu}\eqno
(1.24)$$ where $\mu,\nu \in \Lambda(S(n)^N)$ and
$\xi_{\nu}$ and $\xi_{\mu}$ are prime to one another. Let
$$\Lambda_{\hbox{\sevenrm dom}} = \{\lambda\in \Lambda\mid \lambda\,\,
\hbox{be a
dominant weight}\}. $$

\vs {\bf Remark 3.} By the multiplicity 1-condition note that if
$\nu\in
\Lambda(S(n)^N)$, then $\xi_{\nu}$ is necessarily a homogeneous
polynomial. Define $\hbox{\rm deg}\,\nu$ so that $\xi_{\nu}\in 
 S^{{\hbox{\sevenrm
deg}}\,\nu}(\n)$. Furthermore, clearly $\xi_{\nu}$ is then a
highest weight vector of an irreducible $\g$-module in
$S^{{\hbox{\sevenrm deg}}\,\nu}(\g)$ and in particular
$\nu\in \Lambda_{\hbox{\sevenrm dom}}$. That is, $$\Lambda(S(\n)^N)\s
\Lambda_{{\hbox{\sevenrm dom}}}\cap \Lambda_{\cal B}.\eqno (1.25)$$

\vs {\bf 1.5.} If $\nu\in
\Lambda(S(\n)^N)$, it follows easily from the multiplicity-1 condition 
and the uniqueness of prime factorization that all the prime factors of
$\xi_{\nu}$ are again weight vectors in $S(\n)^N$.  Let
$${\cal P} = \{ \nu \in \Lambda(S(\n)^N)\mid \hbox {$\xi_{\nu}$
be a prime polynomial in $S(\n)^N$}\}.\eqno (1.26)$$ We can then readily prove\vs 

 {\bf Theorem 7.} {\it One has
$\hbox{\rm card}\,{\cal P} = m$ where, we recall $m=\hbox{\rm card}\,{\cal B}$, so that
we can write $${\cal P} = \{\mu_1,\ldots, \mu_{m}\}.\eqno (1.27)$$
Furthermore the weights $\mu_i$ in ${\cal P}$ are
linearly independent and the set $P$ of prime polynomials,
$\xi_{\mu_i},\,i=1,\ldots, m$, are algebraically independent. In
addition, one has a bijection $$\Lambda(S(n)^N) \to (\Bbb
N)^m,\quad\nu \mapsto (d_1(\nu),\ldots, d_m(\nu))\eqno (1.28)$$
such that, writing $d_i= d_i(\nu)$, up to scalar multiplication,
$$\xi_{\nu} = \xi_{\mu_1}^{d_1}\cdots \xi_{\mu_m}^{d_m}\eqno
(1.29)$$ and (1.29) is the prime factorization of $\xi_{\nu}$ for
any $\nu \in \Lambda(S(\n)^N$. Finally, $$S(\n)^{N} = \Bbb
C[\xi_{\mu_1},\ldots,\xi_{\mu_m}]\eqno (1.30)$$ so that $S(\n)^{N}$
is a polynomial ring in $m$-generators.}

\vs {\bf Remark 4.} One may
readily extend part of Theorem 7 to weight vectors in $Q(\n)^N$. In
fact one easily establishes that there is a bijection
$$\Lambda(Q(n_-)^N)\to (\Bbb Z)^m,\,\,\,\gamma\mapsto
(e_1(\gamma),\ldots,e_m(\gamma))$$ so that writing $e_i(\gamma) = e_i$
one has $$\xi_{\gamma} = \xi_{\mu_1}^{e_1}\cdots
\xi_{\mu_m}^{e_m}.\eqno (1.31)$$\vskip .5pc Separating the $e_i$ into 
positive and negative sets yields $\xi_{\nu}$ and $\xi_{\mu}$ of
(1.24).

 \vs {\bf 1.6.} Let $\nu\in \Lambda(S(\n)^N)$. Then by Theorem 6
and (1.25) one has $$\nu\in \Lambda_{\cal B} \cap \Lambda_{\hbox{\sevenrm dom}}$$
so that there exists nonnegative integers $b_{\beta},\,\beta\in {\cal B}$
such that $$\nu = \sum_{\beta\in {\cal B}} b_{\beta}
\beta.\eqno(1.31a)$$ 

\vskip .5pc {\bf Remark 5.} The nonnegativity
follows from dominance since one must have $(\nu,\beta)\geq 0$ for
$\beta\in {\cal B}$.

\vs We wish to prove \vs {\bf Theorem 8.} {\it One has
$$\sum_{\beta\in {\cal B}}\,b_{\beta} = \hbox{\rm deg}\,\nu,\eqno (1.32)$$ and as
a function $\xi_{\nu}\mid \r_-^{\times}$ does not vanish
identically and up to a scalar $$\xi_{\nu}\mid \r_-^{\times}=
\prod_ {\beta\in {\cal B}}\,e_{\beta}^{b_{\beta}}.\eqno (1.33)$$}

\vs
{\bf Proof.} Let
$S^{{\hbox{\sevenrm deg}}\,\nu}(\n)(\nu)$ be the
$\nu$ weight space in $S^{{\hbox{\sevenrm deg}}\,\nu}(\n)$. It does not reduce to
zero since
$\xi_{\nu}\in S^{{\hbox{\sevenrm deg}}\,\nu}(\n)(\nu)$. 
Let $\Gamma$ be the set
of all maps $\gamma:\Delta_+\to \Bbb N$ such that
$$\eqalign{\sum_{\varphi\in
\Delta_{+}}\gamma(\varphi) &= \hbox{\rm deg}\,\nu\cr \sum_{\varphi\in
\Delta_{+}}\gamma(\varphi)\,\varphi &= \nu. \cr}\eqno (1.34)$$
Then if
$$e^{\gamma}= \prod_{\varphi\in \Delta_+}
e_{\varphi}^{\gamma(\varphi)},$$ the set $\{e^{\gamma}\mid \gamma\in
\Gamma\}$ is clearly a basis of $S^{{\hbox{\sevenrm deg}}\,\nu}(\n)(\nu)$ and
consequently unique scalars $s_{\gamma}$ exist so that $$\xi_{\nu} =
\sum_{\gamma\in
\Gamma} s_{\gamma}e^{\gamma}.\eqno (1.35)$$ But by Theorem 5 there
exists
$x\in X$ such that $\xi_{\nu}(x)\neq 0$. However since $X$ is
$B$-homogeneous, the $H$-orbit $\r_-^{\times}$ is contained in $X$ and there
exists $t\in
\r_-^{\times}$ such that $x = u\cdot t$ for some $u\in N$. But since
$\xi_{\nu}$ is $N$-invariant one has $\xi_{\nu}(t)\neq 0$. But from
(1.34) this implies that $$\sum_{\gamma\in
\Gamma} s_{\gamma}e^{\gamma}(t)\neq 0. \eqno(1.36)$$ But
$e^{\gamma}(t)= 0 $ for any $\gamma\in\Gamma$ such that
$\gamma(\varphi)\neq 0$ for $\varphi\notin {\cal B}$. Thus there
exists $\gamma'\in \Gamma$ such that $$\gamma'(\varphi)=0$$ for all 
$\varphi\notin {\cal B}$ and $$e^{\gamma'}(t)\neq 0.\eqno (1.37)$$
But by the independence of ${\cal B}$ one has that $\gamma'$ is
unique and hence one must have $\gamma'(\beta) = b_{\beta}$. A
similar argument yields (1.33). QED
\vs
\centerline{\bf 2. A
representation-theoretic construction, due to Lipsman--Wolf,
}\centerline{\bf of certain elements in $S(\n)^N$}

\vskip 1.5pc {\bf
2.1.} Let $\lambda\in \Lambda_{\hbox{\sevenrm dom}}$ and let $V_{\lambda}$ be
  a
finite-dimensional irreducible $\g$-module with highest
weight$\lambda$. Then, correspondingly, 
$V_{\lambda}$ is a $U(\g)$-module with respect to a surjection
$\pi_{\lambda}:U(\g)\to \hbox{\rm End}\,V_{\lambda}$. Let $0\neq v_{\lambda}\in
V_{\lambda}$ be a highest weight vector. Also let
$V_{\lambda}^*$ be the contragredient dual $\g$-module. The
pairing of $V_{\lambda}$ and $V_{\lambda}^*$ is denoted by
$\langle v,z\rangle$ with $v\in V_{\lambda}$ and $z\in
V_{\lambda}^*$. (We will use this pairing notation throughout in
other contexts.) But as one knows 
$V_{\lambda}^*$ is
$\g$-irreducible with highest weight $\lambda^*\in \Lambda_{\hbox{\sevenrm dom}}$
given by $$\lambda^*= -w_o\lambda. \eqno(2.1)$$ But then by (1.1) and
the mutual orthogonality of roots in the cascade $$-\lambda^* =
\lambda -
\sum_{\beta\in {\cal B}} \lambda(\beta^{\vee})\beta.$$ That is
$$\lambda + \lambda^* = \sum_{\beta\in {\cal
B}}\lambda(\beta^{\vee})\beta\eqno(2.2)$$ and hence
$$\lambda+\lambda^*\in \Lambda_{\cal B}\cap\Lambda_{\hbox{\sevenrm dom}}.\eqno
(2.3)$$ On the other hand, regarding $U(\g)^*$ as a $\g$-module (dualizing the
adjoint action on $U(\g)$) it is clear that if $f\in U(\g)^*$ defined
by putting, for $u\in U(\g)$,  $$f(u) = \langle
u\,v_{\lambda},z_{\lambda^*}\rangle,\eqno (2.4)$$  then $$\eqalign{f\,\,
&\hbox{is $\n$-invariant and}\cr f \,\,&\hbox{is an $\hh$ weight
vector of weight $\lambda + \lambda^*$}.\cr}\eqno (2.5)$$ \vskip .5pc

Now it is true (as will be seen below) that $\lambda + \lambda^*\in
\Lambda(S(\n)^N)$. It is the idea of Lipsman--Wolf to construct
$\xi_{\lambda + \lambda^*}$ using $f$. The method in $[L-W]$ is to
symmetrize $f$ and restrict to $S(\n)$. However Lemma 3.7 in [L-W] is
incorrect (one readily finds counterexamples). But the idea is
correct. One must modify $f$ suitably and this we will do in the next
section.

\vs {\bf 2.2.} Assume $\ss$ is a finite-dimensional Lie
algebra. Let $U_j(\ss),\, j=1,\ldots,$ be the standard filtration of the
enveloping algebra $U(\ss)$. Let $0\neq f\in U(\ss)^*$. We will say
that
$k\geq-1$ is the codegree of $f$ if $k$ is maximal such that $f$
vanishes on $U_{k-1}(\ss)$.  But then if $k$ is the codegree of $f$
and if $x_i\in\ss, \,i=1,\ldots,k$, and $\sigma$ is any permutation
of $\{1,\ldots,k\}$, then $(x_1\cdots x_k -x_{\sigma(1)}\cdots
x_{\sigma(k)})\in U_{k-1}(\ss)$ so that $$f(x_1\cdots x_k) = f(x_{\sigma(1)}\cdots
x_{\sigma(k)}).\eqno (2.6)$$ But this readily implies that there exists 
a unique element $f_{(k)}\in S^k(\ss)$ such that for any $u\in
U_k(\ss)$ one has $$f_{(k)}(\widetilde {u}) = f(u)\eqno (2.7)$$ where
$\widetilde {u}\in S^k(\ss)$ is the image of $u$ under the
Birkhoff--Witt surjection $U_k(\ss)\to S^k(\ss)$. 

Now let $\ss = \g $ and let $f$ be given by (2.4). Let $k$ be the
codegree of $f$. Identify $\g$ with $\g^*$ using the Killing form.
Then $f_{(k)}\in (S^k(g))^N$ and is an $H$-weight vector of weight
$\lambda + \lambda^*$. On the other hand, by (1.2), $$U_k(\g) =
U_k(\n_-)
\oplus U_{k-1}(\g)\b.\eqno (2.8)$$ However $\b\cdot v_{\lambda}\s
\Bbb C\,v_{\lambda}$ so that $ f$ vanishes on $U_{k-1}(\g)\b$. But this
readily implies $f_{(k)} \in S(\n)^N$. We have proved 

\vs {\bf
Theorem 9.} {\it Let $f$ be given by (2.4) and let $k$ be the
codegree of $f$. Then $\lambda +\lambda^*\in \Lambda(S(\n)^N)$.
Furthermore $k = \hbox{\rm deg} (\lambda +\lambda^*)$ and up to scalar
multiplication $$f_{(k)} = \xi_{\lambda + \lambda^*}.\eqno
(2.9)$$}
The inclusion (1.25) is actually an equality
$$\Lambda(S(\n)^N) = \Lambda_{\hbox{\sevenrm dom}}\cap \Lambda_{\cal B}.\eqno
(2.10)$$ This equality is due to Tony Joseph and I was not aware of
it until read it in [J].  However, the equality (2.10) follows
immediately from the modified Lipsman--Wolf construction Theorem 9.
Indeed let $\nu\in \Lambda_{\hbox{\sevenrm dom}}\cap \Lambda_{\cal B}$. To show
$\nu \in \Lambda(S(\n)^N$, it suffices to show that $$e_i(\nu)\geq
0\eqno (2.11)$$ in (1.31) for any $i=1,\ldots, m$. But putting
$\lambda = \nu$, one has $\lambda +\lambda^* = 2\nu$ and by Theorem 9
one has all $e_i(2\nu)\geq 0$. But clearly $e_i(2\nu) = 2e_i(\nu)$.
This proves (2.11). 

\vs
The results in this paper will appear in [K1] in Progress in Mathematics,
 in honor of Joe.

\vskip 2pc

\centerline{\bf References}

\item {[J]} Anthony Joseph, { A preparation theorem for the
prime spectrum of a semisimple Lie algebra}, {\it J. of Alg.}, {\bf
48}, No.2, (1977), 241--289.

\item{[K]} The cascade of orthogonal roots and the coadjoint structure
 of the nilradical of a Borel subgroup of a semisimple Lie group, paper
  in honor of I.M. Gelfand, to appear in the {\it  Moscow Math Journal},
edited by Victor Ginzburg, Spring 2012. 

\item {[K1]} Bertram Kostant, {\it $\hbox{Cent}\, U(\n)$, cascade of orthogonal roots, and a construction
of Lipsman-Wolf}, arXiv /0178653, [math.RT], Jan 12, 2011; to appear in {\it Lie Groups, Structures, 
and  Representation Theory}, in honor of Joseph A. Wolf,  
Progress in Mathematics,   2012.

\item{[LW]} Ronald Lipsman and Joseph A. Wolf, Canonical semi-invariants and
the Plancherel formula for parabolic groups, {\it Trans. Amer. Math. Soc.}
{\bf 269} (1982), 111--131.

\vs
\noindent Bertram Kostant {(Emeritus)}

  \noindent Department of Mathematics

 \noindent  MIT

  \noindent  Cambridge, MA 02139

  \noindent  USA

  \noindent  kostant @math.mit.edu

\end

\end